\documentclass[reqno,11pt]{amsart}
\usepackage{amscd,amssymb,verbatim}

\numberwithin{equation}{section}
\theoremstyle{plain}

\newcommand\alp{\alpha}

         \newcommand\Del{\Delta}
\newcommand\eps{\varepsilon}
\newcommand\zet{\zeta}

\newcommand\lam{\lambda}                \newcommand\Lam{\Lambda}

         \newcommand\Ome{\Omega}

\newcommand\calM{{\mathcal{M}}}

\newcommand\RR{\mathbb{R}}

\newcommand\ZZ{\mathbb{Z}}

\newcommand\CC{\mathbb{C}}

\newcommand\NN{\mathbb{N}}

\newcommand\nek{,\ldots,}
\newcommand\sdp{\times \hskip -0.3em {\raise 0.3ex
\hbox{$\scriptscriptstyle |$}}} 

\newcommand\Id{\operatorname {Id}}

\newcommand\ind{\operatorname{ind}}

\newcommand\Ker{\operatorname{Ker}}

\newcommand\rk{\operatorname{rk}}

\newcommand\RE{\operatorname{Re}}

\newcommand\spec{{\rm spec}}

\newcommand\supp{\operatorname{supp}}

\newcommand\Tr{\operatorname{Tr}}

\newcommand\of{{\overline{f}}}
\newcommand\oF{{\overline{F}}}

\newcommand\om{{\overline{m}}}

\newcommand\oW{{\overline{W}}}

\newcommand\tilA{{\widetilde{A}}}

\newcommand\tilB{{\widetilde{B}}}

\newcommand\tilE{{\widetilde{E}}}
\newcommand\tilf{{\widetilde{f}}}
\newcommand\tilF{{\widetilde{F}}}

\newcommand\tilK{{\widetilde{K}}}

\newcommand\tilM{{\widetilde{M}}}
\newcommand\tilH{{\widetilde{H}}}

\newcommand\tilP{{\widetilde{P}}}

\newcommand\tilU{{\widetilde{U}}}

\newcommand\tilV{{\widetilde{V}}}

\newcommand\tilDel{{\widetilde{\Delta}}}
\newcommand\tilchi{{\widetilde{\chi}}}

\theoremstyle{plain}
\newtheorem{Thm}[subsection]{Theorem}
\newtheorem{Cor}[subsection]{Corollary}
\newtheorem{Lem}[subsection]{Lemma}
\newtheorem{Prop}[subsection]{Proposition}
\newtheorem{Conjec}[subsection]{Conjecture}

\theoremstyle{definition}
\newtheorem{Def}[subsection]{Definition}

\theoremstyle{remark}

\newtheorem{Rem}[subsection]{Remark}

\newif\ifShowLabels
\ShowLabelstrue
\newdimen\theight
\def\TeXref#1{%
        \leavevmode\vadjust{\setbox0=\hbox{{\tt
                \quad\quad  {\small \textrm #1}}}%
        \theight=\ht0
        \advance\theight by \lineskip
        \kern -\theight \vbox to
        \theight{\rightline{\rlap{\box0}}%
        \vss}%
        }}%
\newif\ifShowLabels
\ShowLabelstrue
\newdimen\theight
\def\TeXrefEq#1{%
        \leavevmode\vadjust{\setbox0=\hbox{{\tt
                \quad\quad  {\small \textrm #1}}}%
        \theight=\ht1
        \advance\theight by \lineskip
        \kern -\theight \vbox to
        \theight{\rightline{\rlap{\box0}}%
        \vss}%
        }}%

\newcommand{\refs}[1]{Section ~\ref{S:#1}}
\newcommand{\refss}[1]{Subsection ~\ref{SS:#1}}
\newcommand{\reft}[1]{Theorem ~\ref{T:#1}}
\newcommand{\refl}[1]{Lemma ~\ref{L:#1}}
\newcommand{\refp}[1]{Proposition ~\ref{P:#1}}
\newcommand{\refc}[1]{Corollary ~\ref{C:#1}}

\newcommand{\refd}[1]{Definition ~\ref{D:#1}}

\newcommand{\refe}[1]{\eqref{E:#1}}

\newenvironment{thm}[1]%
        { \begin{Thm} \label{T:#1}  \ifShowLabels \TeXref{T:#1} \fi }%
        { \end{Thm} }

\renewcommand{\th}[1]{\begin{thm}{#1} \sl }
\renewcommand{\eth}{\end{thm} }

\newenvironment{lemma}[1]%
        { \begin{Lem} \label{L:#1}  \ifShowLabels \TeXref{L:#1} \fi }%
        { \end{Lem} }
\newcommand{\lem}[1]{\begin{lemma}{#1} \sl}
\newcommand{\elem}{\end{lemma}}

\newenvironment{propos}[1]%
        { \begin{Prop} \label{P:#1}  \ifShowLabels \TeXref{P:#1} \fi }%
        { \end{Prop} }
\newcommand{\prop}[1]{\begin{propos}{#1}\sl }
\newcommand{\eprop}{\end{propos}}

\newenvironment{corol}[1]%
        { \begin{Cor} \label{C:#1}  \ifShowLabels \TeXref{C:#1} \fi }%
        { \end{Cor} }
\newcommand{\cor}[1]{\begin{corol}{#1} \sl }
\newcommand{\ecor}{\end{corol}}

\newenvironment{conjec}[1]%
        { \begin{Conjec} \label{Conj:#1}  \ifShowLabels \TeXref{C:#1} \fi }%
        { \end{Conjec} }
\newcommand{\conj}[1]{\begin{conjec}{#1} \sl }
\newcommand{\econj}{\end{conjec}}

\newenvironment{defeni}[1]%
        { \begin{Def} \label{D:#1}  \ifShowLabels \TeXref{D:#1} \fi }%
        { \end{Def} }
\newcommand{\defe}[1]{\begin{defeni}{#1} \sl }
\newcommand{\edefe}{\end{defeni}}

\newenvironment{remark}[1]%
        { \begin{Rem} \label{R:#1}  \ifShowLabels \TeXref{R:#1} \fi }%
        { \end{Rem} }
\newcommand{\rem}[1]{\begin{remark}{#1}}
\newcommand{\erem}{\end{remark}}

\newcommand{\eq}[1]%
        { \ifShowLabels \TeXrefEq{E:#1} \fi
           \begin{equation} \label{E:#1} }
\newcommand{\eeq}{\end{equation}}

\newcommand{\prf}{ \begin{proof} }
\newcommand{\eprf}{ \end{proof} }
\newcommand{\Label}[1]{\label{#1}  \ifShowLabels \TeXref{#1} \fi }

\ShowLabelsfalse

\newcommand{\ten}{\otimes}
\renewcommand{\b}{\bullet}
\newcommand{\h}[1]{\text{\( H^{#1}(M,F)\)}}
\newcommand{\hb}{\h{\b}}
\renewcommand{\dh}[1]{\text{\( \det H^{#1}(M,F)\)}}
\newcommand{\dhb}{\text{\( \det H^{\b}(M,F)\)}}

\newcommand{\ts}{\text{\( C^{\b}(W^u,F)\)}}
\newcommand{\dts}{\text{\( \det C^{\b}(W^u,F)\)}}

\renewcommand{\om}{\text{\( \Ome^{\b}(M,F)\)}}

\newcommand{\lapt}{\Del_{f,t}}

\newcommand{\fdm}{\text{\( \|\cdot\|_{\det H^{\b}(V)}\)}}
\newcommand{\fdmb}{\text{\( |\cdot|_{\det H^{\b}(V)}\)}}

\newcommand{\rsm}{\text{\( \|\cdot\|^{RS}_{\det \h{\b}}\)}}

\newcommand{\mm}{\text{\( \|\cdot\|^{\calM}_{\det \h{\b}}\)}}

\newcommand{\rsmb}{\text{\( |\cdot|^{RS}_{\det \h{\b}}\)}}
\newcommand{\rsmbt}{\text{\( |\cdot|^{RS}_{\det \h{\b},f,t}\)}}

\newcommand{\rst}{\text{\( \rho^{RS}\)}}
\newcommand{\rstt}{\text{\( \rho^{RS}(f,t)\)}}
\newcommand{\rstts}{\text{\( \rho^{RS}_{sm}(f,t)\)}}
\newcommand{\rsttl}{\text{\( \rho^{RS}_{la}(f,t)\)}}

\newcommand{\gtm}{\text{\( g^{TM}\)}}
\newcommand{\gft}{\text{\( g^F_t\)}}
\newcommand{\gf}{\text{\( g^F\)}}

\newcommand{\nf}{\text{\( \nabla f\)}}

\newcommand{\p}{\text{\( \partial\)}}

\newcommand{\ld}{\log\det{}'}

\newcommand{\G}{C^\infty(M,E)}

\begin{document}
\title{New proof of the Cheeger-M\"uller Theorem}
\author{Maxim Braverman}
\address{Department of Mathematics\\
        Northeastern University   \\
        Boston, MA 02115 \\
        USA
         }
\email{maxim@neu.edu}
\thanks{Research was partially supported by the NSF grant DMS-0204421}
\begin{abstract}
We present a short analytic proof of the equality between the analytic and combinatorial
torsion. We use the same approach as in the proof given by Burghelea, Friedlander and
Kappeler, but avoid using the difficult Mayer-Vietoris type formula for the determinants
of elliptic operators. Instead, we provide a direct way of analyzing the behaviour of the
determinant of the Witten deformation of the Laplacian. In particular, we show that this
determinant can be written as a sum of two terms, one of which has an asymptotic
expansion with computable coefficients and the other is very simple (no zeta-function
regularization is involved in its definition).
\end{abstract}
\maketitle 
\section{Introduction}\Label{S:setting}

\subsection{Cheeger-M\"uller theorem}\Label{SS:setting}
Let $F$ be a flat vector bundle over a compact odd dimensional Riemannian manifold $M$.
Suppose that $F$ is equipped with a Hermitian metric $g^F$, which induces a flat metric
on the determinant line bundle $\det{F}$. These data define the Ray-Singer metric \rsm\
on the determinant line \dhb, cf. \cite[Def.~2.2]{BisZ1} and \refd{rsmet} of this paper.

Let $f:M\to \RR$ be a Morse function satisfying the Thom-Smale transversality conditions,
\cite{Sm74,Sm73}. Then one defines the Milnor metric \mm\ on \dhb, cf.
\cite[Def.~1.9]{BisZ1} and \refd{milmet} of this paper.
\th{RS}
 \(\displaystyle
    \rsm \ = \ \mm.
 \)
\eth
For the case when the metric $g^F$ is flat, the theorem was conjectured by  Ray and
Singer \cite{RaSi1}. The Ray-Singer conjecture was proven independently by Cheeger
\cite{Ch} and M\"uller \cite{Mu1}. Later M\"uller \cite{Mu2} extended the result to the
case when $g^F$ is not necessarily flat, but the induced metric on $\det{F}$ is flat. The
methods of Cheeger and M\"uller are both based on a combination of the topological and
analytical methods.  Bismut and Zhang \cite{BisZ1} suggested a purely analytical proof of
the Ray-Singer conjecture and generalized it to the case, when the dimension of $M$ is
not necessarily odd and the induced metric on $\det{F}$ is not flat.

Another purely analytical proof of \reft{RS} was suggested by Burghelea, Friedlander and
Kappeler \cite{BFK3}. Their method was based on application of the highly non-trivial
Mayer-Vietoris-type formula for the determinant of an elliptic operator \cite{BFK1}.

In this paper we suggest a new proof of \reft{RS}, which essentially follows the lines of
\cite{BFK3} but is considerably simpler in several steps. In particular, we avoid the use
of Mayer-Vietoris-type formula from \cite{BFK1}.

\subsection{The sketch of the proof}\Label{SS:idpr}
Let $d^F:\Ome^\b(M,F)\to \Ome^{\b+1}(M,F)$ be the de Rham differential. Consider the {\em
Witten deformation} $d^F_t=e^{-tf}d^Fe^{tf}$ of $d^F$ and set
$\Del_{f,t}=d_t^Fd^{F*}_t+d^{F*}_td_t^F$. It was shown by Witten \cite{Witten82} that,
when $t\to\infty$, finitely many eigenvalues of $\Del_{f,t}$ tend to zero (these are, so
called ``small" eigenvalues), while the rest of the eigenvalues tend to infinity (these
are ``large" eigenvalues). The Ray-Singer metric can be expressed, roughly speaking, as
the product of the contribution of the ``small" eigenvalues and the contribution of the
``large" eigenvalues, cf. \refss{defmet}. The proof is naturally divided into the study
of those two contributions.

The contribution of the ``small" eigenvalues is summarized in \reft{BZ}, which was proven
by Bismut and Zhang \cite{BisZ1}. The original proof was based on difficult results of
Helffer and Sj\"ostrand \cite{HelSj4}. Later Bismut and Zhang \cite[\S6]{BisZ2} found a
short and very elegant proof  of this result (see also \cite{BFKM}).

It remains to study the contribution of the ``large" eigenvalues, which we denote by
$\rho^{RS}_{la}(f,t)$, cf. \refss{defmet}. Let $P_{la,t}$ denote the orthogonal
projection on the span of the eigenforms corresponding to the ``large" eigenvalues of
$\Del_{f,t}$. The contribution $\rho^{RS}_{la}(f,t)$ of the ``large" eigenvalues is
defined in terms of $\log\det'[\Del_{f,t}P_{la,t}]$, where $\det'$ denotes the
regularized determinant, cf. \refss{rstor}. Our method is based on the following simple
formula, cf. \refp{aeDel},
\begin{multline}\Label{E:aeDel0}
   k\, \ld\big[\lapt{}P_{la,t}]
   \ = \\
   \ld\big[\, (\Del^k_{f,t}+t^{2k})\, P_{la,t}\big]
   \ - \
    t^{2k}\int_0^1\, \Tr\big[\, (\Del^k_{f,t}+\tau t^{2k})^{-1}P_{la,t}\, \big]\, d\tau.
\end{multline}
Here we choose $k>n/2$, so that the operator  $(\Del^k_{f,t}P_{la,t}+\tau t^{2k})^{-1}$
is of trace class.

The first summand in the right hand side of \refe{aeDel0} is the logarithm of the
determinant of an {\em operator elliptic with parameter}, cf. \cite{Sh1,BFK1}. It is
shown in the Appendix to \cite{BFK1} that it has a nice asymptotic expansion with
computable coefficients. The second summand, though does not have an asymptotic
expansion, is very simple since no $\zet$-function regularization is needed to define it.
It is not difficult now to prove the following result (cf. \reft{compare}): Let $\tilM$
be another Riemannian manifold and $\tilF\to \tilM$ be a flat vector bundle over $\tilM$
such that $\dim{\tilF}= \dim{F}$. Let $\tilf:\tilM\to \RR$ be a Morse function. Assume
that the functions $f$ and $\tilf$ have {\em the same critical points structure}, cf.
\refd{samecp}. Then $\log\rho^{RS}_{la}(f,t)- \log\rho^{RS}_{la}(\tilf,t)$ has a nice
asymptotic expansion with computable coefficients. This result was central in the
Burghelea-Friedlander-Kappeler proof, cf. Theorem~B of \cite{BFK3}.

Set
\eq{R}
    R(M,F,f) \ := \ \log\frac{\rsm}{\mm}.
\end{equation}
By \cite[Th.~9.3]{Milnor1}, the Milnor metric, and, hence, $R(M,F,f)$ is independent of
$f$. It is, however, convenient to keep $f$ in the notation.

In \refs{prRS}, we show, that, if $f$ and $\tilf$ have the same critical points
structure, then $R(M,F,f)= R(\tilM,\tilF,\tilf)$. It follows from \cite{Milnor-hcob} that
there exist Morse functions $f_1, f_2$ satisfying the Thom-Smale condition on
$M\times{S^2}$ and $M\times{S^1}\times{S^1}$ respectively, which have the same critical
points structure. Let $F_1, F_2$ denote the lifts of $F$ to $M\times{S^2}$ and
$M\times{S^1}\times{S^1}$ respectively. Then
\eq{R2}
    R(M\times{S^2},F_1,f_1) \ = \ R(M\times{S^1}\times S^1,F_2,f_2).
\end{equation}
Theorem~2.5 of \cite{RaSi1} expresses the Ray-Singer torsion of the product $M\times{N}$
(here $N$ is a compact Riemannian manifold) in terms of the Ray-Singer torsion of $M$. In
\refs{prRS}, we use this result to show that
\eq{1-2}
    R(M,F,f) \ = \ R(M\times{S^2},F_1,f_1), \qquad\text{and}\qquad
    R(M\times{S^1}\times S^1,F_2,f_2)  \ = \ 0.
\end{equation}
Combining \refe{R2} and \refe{1-2} we obtain $R(M,F,f)=0$.

\subsection{The results used in the proof}\Label{SS:used}
For convenience of the reader, we list all the results which we use but don't prove in
this paper.
\begin{itemize}
\item Topological invariance of the Milnor and the Ray-Singer torsion. The proofs can be
found in \cite[Th.~9.3]{Milnor1} and \cite[Th.~2.1]{RaSi1} respectively.

\item The relationship between the Milnor metric and the contribution of the ``small"
eigenvalues of the Witten deformation of the Laplacian to the Ray-Singer torsion, cf.
\reft{BZ}. A very nice proof can be found in \cite[\S6]{BisZ2} (see also \cite{BFKM}).

\item The asymptotic expansion of the trace and the determinant of an operator elliptic
with parameter obtained in the Appendix to \cite{BFK1}.

\item Existence of a constant $C>0$ such that
$\Tr\big[(\Del^{k}_{f,t}+\eps)^{-1}P_{la,t}\big]< C$ for all $k>n$, $\eps>0$ and
$|t|\gg0$. This simple estimate follows, for example, from Lemma~3.3 of \cite{BFK3}.

\item The expression for Ray-Singer torsion on the product of 2 manifolds, cf.
\cite[Th.~2.5]{RaSi1}.
\end{itemize}
Apart from these results the paper is completely independent.

\subsection*{Acknowledgements}
I would like to thank Alexander Abanov for valuable discussions.

\section{The Milnor metric and the Milnor torsion}\Label{S:milmet}

%
\subsection{The determinant line of a finite dimensional complex}\Label{SS:detline}
If $\lam$ is a real line, let $\lam^{-1}$ be the dual line. If $E$ is a finite
dimensional vector space, set
 \(
    \det E= \Lam^{\max}(E).
 \)
Let \/
 \(
    (V^\b,\p):\ 0\to V^0\to\cdots\to V^n\to 0
 \) \/
be a  complex of finite dimensional Euclidean vector spaces. Let
$H^\b(V)=\bigoplus_{i=0}^n H^{i}(V)$ be the cohomology of $(V^{\b},\p)$. Set
\[
    \det V^\b=\bigotimes_{i=0}^n \big(\det V^{i}\big)^{(-1)^i},
    \qquad
    \det H^{\b}(V)=\bigotimes_{i=0}^n \Big(\det H^{i}(V)\Big)^{(-1)^i}.
\]
Then, by \cite{KnMu}, there is a canonical isomorphism of real lines
\eq{kmu-isom}
    \det H^{\b}(V)\simeq \det V^\b.
\end{equation}

\subsection{Two metrics on the determinant line}\Label{SS:2met}
The Euclidean structure on $V^\b$ defines a metric on $\det V^\b$. Let $\fdm$ be the
metric on the line $\det H^\b(V)$ corresponding to this metric via the canonical
isomorphism \refe{kmu-isom}.

Let $\p^*$ be the adjoint of $\p$ with respect to the Euclidean structure on
$V^\b$. Using the finite dimensional Hodge theory, we have the canonical
identification
\eq{hodge-ts}
   H^i(V^\b,\p)\simeq \{v\in V^i:\ \p v=0, \p^* v=0\},\ \ 0\le i\le n.
\end{equation}
As a vector subspace of $V^i$, the vector space in the right-hand side of \refe{hodge-ts}
inherits the Euclidean metric. We denote by \fdmb\ the corresponding metric on $\det
H^{\b}(V)$.

The metrics \fdm\ and \fdmb\ do not coincide in general. We shall describe the
discrepancy.

Set $\Del=\p\p^*+\p^*\p$ and let $\Del^i$ denote the restriction of $\Del$ to $V^i$. Let
\/ $\det'\Del^i$ \/ denote the product of the non-zero eigenvalues of $\Del^i$.
\defe{fdtor}
The {\em  torsion} $\rho$ of the complex $(V^\b,d)$ is the number defined by the formula
\[
    \log\rho \ = \ \frac12\sum_{i=0}^n\, (-1)^ii\log\det{}'\Del^i.
\]
\edefe
The following result is proved, e.g., in \cite[Prop.~1.5]{BisGS}
\[
    \fdm=\fdmb \cdot\rho.
\]

\subsection{The Thom-Smale complex}\Label{SS:th-sm}
Let $f:M\to\RR$ be a Morse function satisfying the Smale transversality conditions
\cite{Sm74,Sm73}\ (for any two critical points $x$ and $y$ of $f$ the stable manifold
$W^s(x)$ and the unstable manifold $W^u(y)$, with respect to \nf, intersect
transversely).

Let $B$ be the set of critical points of $f$.   If $x\in B$, let $F_x$ denote the fiber
of $F$ over $x$ and let $[W^u(x)]$ denote the real line generated by $W^u(x)$. For $0\le
i\le n$, set
\[
  C^i(W^u,F)=
    \bigoplus_{\begin{smallmatrix} x\in B\\ \ind(x)=i
    \end{smallmatrix}}
        [W^u(x)]^*\ten_{\RR} F_x.
\]
By a basic result of Thom (\cite{Th1}) and Smale (\cite{Sm73}) (see also
\cite[pp.~28--30]{BisZ1}), there are well defined linear operators
\[
    \p:C^i(W^u,F)\to C^{i+1}(W^u,F),
\]
such that  the pair $(C^\b(W^u,F),\p)$ is a complex  and there is a canonical
identification of $\ZZ$-graded vector spaces
\eq{tsc=c}
      H^\b(C^\b(W^u,F),\p)\simeq H^\b(M,F).
\end{equation}
\subsection{The Milnor metric}\Label{SS:milmet}
By \refe{kmu-isom} and \refe{tsc=c}, we know that
\eq{km-ts}
    \dh{\b}\simeq \dts.
\end{equation}
The metric \gf\ on $F$ determines the structure of an Euclidean vector space on \ts. This
structure induces a metric on $\det\ts$.
\defe{milmet}
   The {\em Milnor metric} \mm\ on the line \dh{\b}\ (cf.
   \cite[\S{}Id]{BisZ1}) is the metric
   corresponding to the above metric on \dts\ via the
   canonical isomorphism \refe{km-ts}.
\edefe
By Milnor \cite[Th.~9.3]{Milnor1},  the Milnor metric coincides with the Reidemeister
metric defined through a smooth triangulation of $M$. It follows that  \mm\ does not
depend upon $f$ and \gtm, \gf\ and, hence, is a topological invariant of the flat bundle
$F$.

\section{The Ray-Singer metric and the Ray-Singer torsion}\Label{S:rsmet}

\subsection{The $L^2$ metric on the determinant line}\Label{SS:badrsmet}
Let $(\om,d^F)$ be the de Rham complex of the smooth sections of $\Lam(T^*M)\ten F$
equipped with the coboundary operator $d^F$. The cohomology of this complex is
canonically isomorphic to \hb.

Let $\ast$ be the Hodge operator associated to the metric \gtm. We equip \om\ with the
inner product
\eq{dr-inpr}
    \langle\alp,\alp'\rangle_{\om}=
        \int_M\langle\alp\land\ast \alp'\rangle_{\gf}.
\end{equation}
By Hodge theory, we can identify \hb\ with the space of harmonic forms in \om. This space
inherits the Euclidean product \refe{dr-inpr}. The {\em $L^2$-metric} \rsmb\ on \dhb\ is
the metric induced by this product.

\subsection{The Ray-Singer torsion}\Label{SS:rstor}
Let $d^{F*}$ be the formal adjoint of $d^F$ with respect to the metrics \gtm\ and \gf.
Let $\Del=d^Fd^{F\ast}+d^{F\ast}d^F$ be the Laplacian and let $\Del^i$ denote the
restriction of $\Del$ to $\Ome^i(M,F)$. Let $P^i:\om\to\Ker\Del^i$ be the orthogonal
projection.

To define the torsion of the complex $(\om,d^F)$ one needs to make sense of the notion of
determinant of the Laplacian. This is done using the zeta-function regularization as
follows.

{}For $s\in\CC$, $\RE\, s >n/2$, set
 \(
    \zet^{RS}_{i}(s)= -\Tr\big[(\Del^i)^{-s}(I-P^i)\big].
 \)
By a result of Seeley \cite{See1}, $\zet^{RS}_{i}(s)$ extends to a meromorphic function
of $s\in\CC$, which is holomorphic at $s=0$. Define the determinant $\det'\Del^i$ by the
formula
\[
    \ld\Del^i \ = \ \frac{d}{ds}\zet^{RS}_{i}(0).
\]

\defe{rsmet}
The {\em Ray-Singer torsion \/} $\rst$ is defined by the formula (cf.
\cite[Def.~2.2]{BisZ1})
\[
    \log\rst \ = \ \frac12\sum_{i=0}^n\, (-1)^ii\ld\Del^i.
\]
The {\em Ray-Singer metric \/} \rsm\ on the line \dhb\ is the product
\[
  \rsm=\rsmb\cdot\rst.
\]
\edefe
Ray and Singer \cite[Th.~2.1]{RaSi1} proved that the metric \rsm\ is a topological
invariant, i.e., does not depend on the metrics \gtm\ or \gf.

\section{The Witten deformation.}\Label{S:witten}

\subsection{A simplifying assumption}\Label{SS:ass}
Recall that $g^{TM}$ denotes the Riemannian metric on $M$. Following \cite{BFK3} we give
the following
\defe{crpoint}
The pair $(g^{TM},f)$ is a {\em generalized triangulation} of $M$, if $f$ is a
Morse function $f:M\to \RR$ satisfying the Thom-Smale transversality condition
(cf. \refss{th-sm}) and in a neighborhood of every critical point $x$ of $f$
one can introduce local coordinates $(y_1\nek y_n)$ such that
\[
    f(y) \ = \ f(x)  \ - \ \frac12(y_1^2+ \cdots + y_k^2) \ + \
        \frac12(y_{k+1}^2+ \cdots + y_n^2),
\]
and the metric $g^{TM}$ is Euclidean in these coordinates.
\edefe
Since both the Milnor and the Ray-Singer metrics are independent of the choice of $f$ and
$g^{TM}$, it is enough to prove \reft{RS} for the case when $(g^{TM},f)$ is a generalized
triangulation, which we will henceforth assume.

\subsection{The Witten deformation of the Laplacian}\Label{SS:witdef}
Set $d^F_t=e^{-tf}d^Fe^{tf}$, \ $d^{F*}_t=e^{tf}d^Fe^{-tf}$. Then $d^{F*}_t$ is the
formal adjoint of $d^F_t$ with respect to the scalar product \refe{dr-inpr}. The {\em
Witten Laplacian} is the operator
\[
    \lapt \ = \ d^F_td^{F*}_t \ + \ d^{F*}_td^F_t.
\]
We denote by $\lapt^i$ the restriction of $\lapt$ to $\Ome^i(M,F)$.
Let $\rstt$ be the torsion defined as in \refss{rstor}, but with replacing everywhere
$\Del$ by $\lapt$.

The following theorem is well known, cf. \cite{Witten82}
\th{split}
Suppose that the pair $(g^{TM},f)$ is a generalized triangulation.

1. \ There exist positive constants $C',C''$, and $t_0>1/C''$, so that for $|t|\ge t_0$,
we have $\spec(\lapt)\subset [0,e^{-|t|C'})\cup(C''|t|,\infty)$.

2. \ Let $E_{sm,t}^\b\subset \Ome^\b(M,F)$ denote the span
of the eigenvectors of $\lapt$ with eigenvalues less than
$e^{-|t|C'}$. Then $\dim{}E^i_{sm,t}=m_i\rk(F), \ i=1\nek
n$, where $m_i$ is the number of the critical points of $f$
with index $i$.
\eth
Clearly, $E_{sm,t}^\b$ is a subcomplex of the complex $(\Ome^\b(M,F),d^F_t)$. Let
$\rstts$ be the torsion of this subcomplex.  Let $P_{sm,t}^i:\Ome^i(M,F)\to E^i_{sm,t}$
be the orthogonal projection and let $P_{la,t}^i=1-P_{sm,t}^i$. Set
$\zet^{RS}_{la,i}(s)=-\Tr\big[(\lapt^i)^{-s}P_{la,t}^i\big]$ and
\eq{rsttl}
 \begin{aligned}
      \ld\big[\lapt^iP_{la,t}^i] \ &= \ \frac{d}{ds}\zet^{RS}_{la,i}(0); \\
    \log\rsttl \ = \ \frac12\sum_{i=0}^n\, &(-1)^ii\ld\big[\lapt^iP_{la,t}^i].
 \end{aligned}
\end{equation}
Clearly,
\eq{r-smla}
    \rstt \ = \ \rsttl\cdot\rstts \qquad\text{for}\quad
    |t|>t_0.
\end{equation}

\subsection{The Witten Laplacian and the Ray-Singer metric}\Label{SS:defmet}
For each $t\in\RR$, consider the metric $\gft=e^{-2tf}\gf$. Let \rsmbt\  be the
$L^2$-metric on \dhb\ associated to the metrics \gft\ and \gtm. The Laplacian
$\tilDel_{f,t}$ associated to the metrics \gft\ and \gtm\ is conjugate to $\lapt$. More
precisely, we have, $\tilDel_{f,t}=e^{tf}\lapt{}e^{-tf}$, cf. \cite[Prop.~5.4]{BisZ1}.
Hence, $\rstt$ equals the Ray-Singer torsion associated to the  metrics \gft\ and \gtm.
Since the Ray-Singer metric is a topological invariant of $F$, it follows that
\eq{rsm=}
    \rsm \ = \ \rsmbt\cdot \rstt, \qquad\text{for any}\quad
    t\in\RR.
\end{equation}

\th{BZ}
Suppose that the pair $(g^{TM},f)$ is a generalized triangulation, cf. \refd{crpoint}.
Then, as $t\to+\infty$, we have
\eq{BZ}
    \log\, \frac{\rsmb\rstts}{\mm} \ = \
        -t\rk(F)\Tr_s^B[f]+\frac12\tilchi'(F)\log\, \Big(\,
            \frac{t}{\pi}\, \Big) \ + \ o(1),
\end{equation}
where $\Tr^B_s[f]=\sum_{x\in B}(-1)^{\ind(x)}f(x)$ and $\tilchi'(F)= \rk(F)\sum_{x\in
B}\, (-1)^{\ind(x)}\ind(x)$.
\eth
The theorem was first proven in \cite[Th.~7.6]{BisZ1} using the difficult results of
Helffer and Sj\"ostrand \cite{HelSj4}. A short and very elegant proof was found by Bismut
and Zhang \cite[\S6]{BisZ2} (see also \cite{BFKM}).

Recall that the number $R(M,F,f)$ was defined in \refe{R}. Using \refe{r-smla} and
\refe{rsm=}, we obtain the following corollary of \reft{BZ}:
\cor{RS=}
There exists a constant $R=R(M,F,f)$, such that
\eq{RS=}
    \log\rsttl \ = \  R(M,F,f)+t\rk(F)\Tr_s^B[f]-\frac12\tilchi'(F)\log\, \Big(\,
            \frac{t}{\pi}\, \Big) \ + \ o(1),
    \qquad t\to+\infty.
\end{equation}
$R$ is independent of $f$ and \reft{RS} is equivalent to the equality $R=0$.
\ecor
Thus the proof of \reft{RS} is reduced to the study of the asymptotic expansion of
$\rsttl$ as $t\to\infty$.

\section{The comparison theorem}\Label{S:almas}



\subsection{}\Label{SS:2mflds}
Let $M,\tilM$ be Riemannian manifolds of the same odd dimension $n$. Let $F, \tilF$ be
flat vector bundles over $M$ and $\tilM$ respectively, equipped with Hermitian metrics,
such that the induced metrics on $\det{}F$ and $\det{}\tilF$ are flat. We assume that
$\dim{}F= \dim\tilF$.

\defe{samecp}
The Morse functions $f:M\to\RR$ and $\tilf:\tilM\to\RR$ have {\em the same critical
points structure} if there exist open neighborhoods $U\subset M, \ \tilU\subset\tilM$ of
the sets of critical points of $f, \tilf$ respectively, and an isometry $\phi:U\to
\tilU$, such that $f=\tilf\circ\phi$.
\edefe
\defe{asexp}
We say that a function $l(t)$ has a nice asymptotic expansion as $t\to\pm\infty$ if
\[
    l(t) \ = \ \sum_{j=0}^n\, a_j(t/|t|)t^j \ + \
            \sum_{k=0}^n\, b_j(t/|t|)t^j\log|t| \ + \ o(1),
\]
and the coefficient $a_0$ (the free term) satisfy the equality $a_0(1)+ a_0(-1)=0$.
\edefe
The main result of this section is the following
\th{compare}
Let $f:M\to\RR, \tilf:\tilM\to\RR$ be Morse functions with the same critical points
structure and let $U,\tilU$ be as in \refd{samecp}. Then the difference $\log\rsttl-
\log{}\rho^{RS}_{la}(\tilf,t)$ has a nice asymptotic expansion.
\eth
The rest of this section is devoted to the proof of \reft{compare}.

\subsection{Determinant of an operator almost elliptic with parameter}\Label{SS:det-ae}
It is more convenient to work in a slightly more general
situation. Suppose $E$ is a Hermitian vector bundle over a
compact Riemannian manifold $M$ of dimension $n$. Consider
the operator
\eq{Ht}
    H_t \ := \ A+tB+t^2V:\, \G \ \to \ \G, \qquad t\in\RR,
\end{equation}
where $A:\G\to\G$ is a second order self-adjoint elliptic differential operator with
positive definite leading symbol, $B=B(x), V=V(x):E\to E$ are self-adjoint bundle maps
and $V(x)\ge0$ for all $x\in M$. Suppose that there exist constants $t_0, c_1, c_2>0$
such that for all $|t|>t_0$, there are finitely many eigenvalues of $H_t$, which are
smaller than $e^{-c_1|t|}$, while all the other eigenvalues of $H_t$ are larger than
$c_2|t|$. Let $P_t$ denote the orthogonal projection onto the span of the eigensections
of $H_t$ with eigenvalues greater than $1$.

Note that $\rk(\Id-P_t)$ is equal to the number of eigenvalues of $H_t$ (counting
multiplicities) which are smaller than $1$. Hence, the function $t\mapsto \rk(\Id-P_t)$
is locally constant for $|t|> \max\{t_0,1/c_2\}$. Set
\eq{m}
    m_\pm \ := \ \rk (\Id-P_t), \qquad \pm t>\max\{t_0,1/c_2\}.
\end{equation}

Assume, in addition, that there exist constants $k>n/2$ and $C>0$ such that
\eq{Tr<C}
    \Tr\, \big[\, (H_t^k+\eps)^{-1}P_t\, \big] \ < \ C,
    \qquad\text{for all}\quad \eps>0, |t|\gg0.
\end{equation}
Note that when $H_t=\Del_{t,f}$ this assumption is satisfied  for every $k>n$ by
\cite[Lemma~3.3]{BFK3}.

We are interested in the behaviour of the function
 \(
    l(t)=  \ld H_tP_t,
 \)
as $t\to\pm\infty$. Note that, if $V(x)>0$ for all $x\in M$, then $H_t$ is an {\em
elliptic operator with parameter}, cf. \cite[Ch.~1]{Sh1}, \cite[Appendix]{BFK1}. Then
$l(t)$ has a nice asymptotic expansion as $t\to\infty$ with computable coefficients, cf.
\cite[Appendix]{BFK1}. If $V(x)$ is not strictly positive for some $x\in M$, this
asymptotic expansion does not hold any more. However, the following result is true: Let
$\tilE$ be a Hermitian vector bundle over another compact Riemannian manifold $\tilM$. We
assume that the rank of $E$ is equal to the rank of $\tilE$. Let
\[
    \tilH_t \ = \
     \tilA+ t\tilB + t^2\tilV: C^\infty(\tilM,\tilE)\to C^\infty(\tilM,\tilE)
\]
be as above. Then $1$ is not an eigenvalue of $\tilH_t$ for $|t|\gg0$. Let $\tilP_t$ be
the orthogonal projector  onto the span of eigensections of $H_t$ with eigenvalues
greater than $1$.
\th{compare2}
Suppose there exist open sets $U\subset M$ and $\tilU\subset \tilM$ such that $V(x)>0$
for all $x\in M\backslash{U}$ and $\tilV(x)>0$ for all $x\in \tilM\backslash{\tilU}$. Let
$\phi:U\to \tilU$ be a diffeomorphism which preserves the Riemannian metric. Assume that
$\psi:\phi^*\tilE|_{\tilU}\to E|_U$ is an isometry, which identifies the restrictions of
$H_t$ to $U$ and of $\tilH_t$ to $\tilU$. Then the function $\ld{}H_tP_t-
\ld{}\tilH_t\tilP_t$ has a nice asymptotic expansion.
\eth
Clearly, \reft{compare} is an immediate consequence of \reft{compare2}. We pass now to
the proof of \reft{compare2}. First we establish the following
\prop{aeDel}
For every $k> n/2$, the following equality holds
\begin{multline}\Label{E:aeDel}
   k\, \ld H_tP_t \ = \ \ld H_t^kP_t
   \ = \\
   \ld\big[\, (H_t^k+t^{2k})\, P_{t}\big]
   \ - \
    t^{2k}\int_0^1\, \Tr\big[\, (H_t^k+\tau t^{2k})^{-1}P_{t}\, \big]\, d\tau.
\end{multline}
\eprop
\prf
For $k> n/2$ the operator $\big[(H_t^k+\tau t^{2k})P_{t}\big]^{-1}$ is of trace class.
Hence
\begin{multline}\notag
    \frac{d}{d\tau}\ld\big[\, (H_t^k+t^{2k})\, P_{t}\big]
    \\ = \
    \Tr\big[\, (H^k_t+\tau t^{2k})^{-1}
            \frac{d}{d\tau}(H_t^k+\tau t^{2k})P_{t}\, \big]
    \ = \
    t^{2k}\Tr\big[\, (H_t^k+\tau t^{2k})^{-1}P_t\, \big].
\end{multline}
Integrating this equality, we obtain \refe{aeDel}.
\eprf

From now on we assume that $k$ is as in \refe{Tr<C}. Then using the definition of $P_t$
we get
\eq{Tr=Tr+1}
    \int_0^1\, \Tr\big[\, (H_t^k+\tau t^{2k})^{-1}P_{t}\, \big]\, d\tau
    \ = \
     \int_{0}^1\,
       \Tr\big[\, (H_t^k+\tau t^{2k}+|t|^{-k})^{-1}P_{t}\, \big]\, d\tau \ + \ o(t^{-2k}),
\end{equation}
as $t\to\infty$.

Recall that the numbers $m_\pm$ were defined in \refe{m}. Clearly
\eq{ldDP}
\begin{aligned}
    \int_{0}^1\,
       \Tr\big[\, (H_t^k&+\tau t^{2k}+|t|^{-k})^{-1}P_{t}\, \big]\, d\tau
     \\ &= \
     \int_{0}^1\,
       \Tr\, (H_t^k+\tau t^{2k}+|t|^{-k})^{-1}\, d\tau
     \ - \ 3km_\pm t^{-2k}\log|t|
     \ + \ o(t^{-2k}),  \\
    \ld\big[\, &(H_t^k+t^{2k})P_{t}\, \big]
     \ = \
     \ld\big(\, H_t^k+t^{2k}\, \big)
        \ - \
                 2km_\pm\log |t| \ + \ o(1),
\end{aligned}
\end{equation}
as $t\to\pm\infty$.

It is shown in the Appendix to \cite{BFK3} that $ \ld{}(H_t^k+t^{2k})$ has a nice
asymptotic expansion. Hence, \reft{compare2} follows from \refe{aeDel}, \refe{Tr=Tr+1},
\refe{ldDP} and the following
\prop{compare3}
Under the assumptions of \reft{compare2} the function
\eq{compare3}
    t^{2k}\, \int_0^1\, \big[\, \Tr\,(H_t^k+\tau{}t^{2k}+|t|^{-k})^{-1} \ - \
        \Tr\,(\tilH_t^k+\tau{}t^{2k}+|t|^{-k})^{-1}\, \big]\, d\tau
\end{equation}
has a nice asymptotic expansion.
\eprop
The rest of this section is occupied with the proof of \refp{compare3}.

\subsection{Notations}\Label{SS:auxiliary}
Let $W\subset M$ be an open set whose closure $\oW\subset U$ and such that $V(x)>0$ for
all $x\not\in W$. Fix $\eps>0$ such that $V(x)>\eps$ for all $x\not\in W$. We can and we
will assume that $\eps<1$. Let $v:M\to [0,\eps]$ be a smooth function such that
$\supp{}v\subset U$ and $v|_W\equiv\eps$. Set
\eq{A}
    A_{t,\tau} \ := \ H_t^k+\tau t^{2k}+|t|^{-k},
    \qquad
    A_{t,\tau,v}  \ := \ H_t^k+\tau t^{2k}+|t|^{-k}+v^2t^{2k}.
\end{equation}

To simplify the notation we will identify $U$ and $\tilU$ via the diffeomorphism
$\phi:U\to \tilU$. In particular, we will consider $v$ as a function on $\tilM$. We
define operators $\tilA_{t,\tau}$ and $\tilA_{t,\tau,v}$ as in \refe{A} but using
$\tilH_t$ instead of $H_t$.

\lem{asexp+v}
Let $K_{\tau,v}(t,x,y)$ denote the Schwartz kernel of the operator $A_{t,\tau,v}^{-1}$.
Then for each $N\in \NN$
\eq{asexp+v}
    K_{\tau,v}(t,x,x)
    \ = \
    \sum_{j=0}^N\, \alp_j(\tau,t/|t|,x) t^{n-j-2k} \ + \ r_{N,\tau}(t,x),
\end{equation}
where $t^Nr_{N,\tau}(t,x)\to 0$ as $t\to\pm\infty$ uniformly in $\tau\in [0,1], \, x\in
M$. The coefficients $\alp_j(\tau,\pm1,x)$ depend continuously on $\tau\in[0,1]$ and can
be expressed in terms of the full symbol of $H_t$ and a finite number of its derivatives.
If $j=2i$ is even, then $\alp_j(\tau,1,x)+ \alp_j(\tau,-1,x)= 0$.
\elem
\prf
Clearly, $A_{t,\tau,v}= H_t^k+\tau t^{2k}+|t|^{-k}+v^2t^{2k}$ is an operator elliptic
with parameter, cf. \cite[Ch.~1]{Sh1}, \cite[Appendix]{BFK1}. The lemma is a consequence
of the standard construction of the parametrix of an operator elliptic with parameter. It
follows immediately, for example, from Lemma~A.8 in \cite{BFK1}.
\eprf
Since $\Tr\, A_{t,\tau,v}^{-1}= \int_M K_{\tau,v}(t,x,x)dx$ we obtain the following
\cor{asexp+v}
The function $t^{2k}\,\int_{0}^1\,\Tr\, A_{t,\tau,v}^{-1}\, d\tau$ has a nice asymptotic
expansion.
\ecor

\refp{compare3} (and, hence, Theorems~\ref{T:compare2} and \ref{T:compare}) follows now
from the following
\lem{compare4}
Under the assumptions of \reft{compare2} we have
\begin{equation}\Label{E:compare4}
      \Tr\, \big[\, A_{t,\tau}^{-1} - A_{t,\tau,v}^{-1}\, \big]
      \ - \
      \Tr\, \big[\, \tilA_{t,\tau}^{-1} - \tilA_{t,\tau,v}^{-1}\, \big]
     \ = \
     o(t^{-2k})
\end{equation}
as $t\to \infty$ uniformly in $\tau\in [0,1]$.
\elem
\prf
We have
\eq{A-Av}
    A_{t,\tau}^{-1} - A_{t,\tau,v}^{-1}
    \ = \
    A_{t,\tau}^{-1}\,  v^2t^{2k}\, A_{t,\tau,v}^{-1}
    \ = \
    A_{t,\tau,v}^{-1}\,  v^2t^{2k}\, A_{t,\tau}^{-1}.
\end{equation}
Hence
\begin{multline}\Label{E:TrA-Av}
      \Tr\, \Big[\, A_{t,\tau}^{-1} - A_{t,\tau,v}^{-1}\, \Big]
      \ = \
      \Tr\, \Big[\, A_{t,\tau}^{-1} \, v^2t^{2k}\, A_{t,\tau,v}\, \Big]
      \ = \
       t^{2k}\, \Tr\, \Big[\, v\, A_{t,\tau,v}^{-1} A_{t,\tau}^{-1} v\, \Big]
       \\ = \
        t^{2k}\, \Tr\, \Big[\, v\, A_{t,\tau,v}^{-2} v\, \Big]
       \ + \
        t^{4k}\, \Tr\, \Big[\, v\, A_{t,\tau,v}^{-2}\, v^2\,
         A_{t,\tau}^{-1} v\, \Big].
\end{multline}
Similar equality is true for $\Tr\Big[\tilA_{t,\tau}^{-1}- \tilA_{t,\tau,v}^{-1}\Big]$.

Let $\tilK_{\tau,v}(t,x,y)$ denote the Schwartz kernel of the operator
$\tilA_{t,\tau,v}^{-1}$. By \refl{asexp+v}, for all $x\in \supp{}v\subset U$ and all
$N\in\NN$ we have $K_{\tau,v}(t,x,x)- \tilK_{\tau,v}(t,x,x)= o(t^{-N})$ as $t\to\infty$
uniformly in $\tau\in[0,1]$. Hence,
\eq{Av-tilAv}
    \Tr\, \Big[\, v\, A_{t,\tau,v}^{-2} v\, \Big]
    \ - \
    \Tr\, \Big[\, v\, \tilA_{t,\tau,v}^{-2} v\, \Big]
    \ = \
    \int_M\, v\, \big(\, K_{\tau,v}(t,x,x)-\tilK_{\tau,v}(t,x,x)\, \big)\, v\, dx
    \ = \ o(t^{-N}),
\end{equation}
as $t\to\infty$ uniformly in $\tau\in[0,1]$.

Let $I_{t,\tau,v}$ denote the left hand side of \refe{compare4}. From \refe{TrA-Av} and
\refe{Av-tilAv} we conclude that
\eq{Ittauv}
    I_{t,\tau,v}
    \ := \
   t^{4k}\, \Tr\, \Big[\, v\, A_{t,\tau,v}^{-2}\, v^2\, A_{t,\tau}^{-1} v\, \Big] \ - \
   t^{4k}\, \Tr\, \Big[\, v\, \tilA_{t,\tau,v}^{-2}\, v^2\, \tilA_{t,\tau}^{-1} v\, \Big]
    \ + \ o(t^{2k-N})
\end{equation}
as $t\to\infty$ uniformly in $\tau\in[0,1]$.

Using the isometry $\psi:\phi^*\tilE|_{\tilU}\to E|_U$ we can view
$v\tilA_{t,\tau,v}^{-2}v$ and $v\tilA_{t,\tau}^{-2}v$ as operators  acting on the space
of sections of the bundle $E$. Then
\begin{multline}\Label{E:separate}
    \left|\,
    \Tr\, \Big[\, v\, A_{t,\tau,v}^{-2}\, v^2\, A_{t,\tau}^{-1} v\, \Big] \ - \
    \Tr\, \Big[\, v\, \tilA_{t,\tau,v}^{-2}\, v^2\, \tilA_{t,\tau}^{-1} v\, \Big]
    \, \right|
    \\ \le \
     \Big|\,
    \Tr\, \Big[\, \big(\,
                     v\, A_{t,\tau,v}^{-2}\, v - v\, \tilA_{t,\tau,v}^{-2}\, v\,
                   \big)\, v\,  A_{t,\tau}^{-1} v\, \Big] \, \Big|
    \ + \
    \Big|\,
     \Tr\, \Big[\, v\, \tilA_{t,\tau,v}^{-2}\, v \,
       \big(\, v\, A_{t,\tau}^{-1}\, v - v\, \tilA_{t,\tau}^{-1}\, v\, \big)
       \Big] \, \Big|
     \\ \le \
         \big\|\, v\,  A_{t,\tau}^{-1} v\, \big\|\cdot \Big|\, \Tr\, \Big[\,
       v\, A_{t,\tau,v}^{-2}\, v - v\, \tilA_{t,\tau,v}^{-2}\, v\, \Big] \, \Big|
    \ + \
    \big\|\, v\, \tilA_{t,\tau,v}^{-2}\, v \, \big\|\cdot
    \Big|\, \Tr\, \Big[\, v\, A_{t,\tau}^{-1}\, v -
             v\, \tilA_{t,\tau}^{-1}\, v\, \Big]\, \Big|.
\end{multline}

Fix $N>7k$. Then, using \refe{Av-tilAv}, \refe{Ittauv}, \refe{separate} and the obvious
estimates
\eq{norms>}\notag
    \big\|\, v\, A_{t,\tau}^{-1}\, v\, \big\| \ \le \ |t|^k,
    \qquad
    \big\|\, v\, \tilA_{t,\tau,v}^{-2}\, v\, \big\|
     \ \le \ \eps^{-2}\, t^{-4k},
\end{equation}
we conclude
\begin{equation} \Label{E:I=}
    \big|\, I_{t,\tau,v}\, \big| \ \le  \
       \eps^{-2}\, \Big|\, \Tr\, \Big[\, v\, A_{t,\tau}^{-1}\, v -
             v\, \tilA_{t,\tau}^{-1}\, v\, \Big]\, \Big| \ + \ o(t^{-2k}).
\end{equation}

Applying again \refe{A-Av}, we get
\begin{multline}\Label{E:vA-Avv}
    v\, A_{t,\tau}^{-1}\, v - v\, \tilA_{t,\tau}^{-1}\, v
    \ = \
   t^{2k}\, \Big(\, v\, A_{t,\tau,v}^{-1}\, v^2\, A_{t,\tau}^{-1}\, v -
         v\, \tilA_{t,\tau,v}^{-1}\, v^2\, \tilA_{t,\tau}^{-1}\, v\, \Big)
    \ + \
    \Big(\, v\, A_{t,\tau,v}^{-1}\, v -
             v\, \tilA_{t,\tau,v}^{-1}\, v\, \Big)
   \\ = \
   t^{2k}\, v\, A_{t,\tau,v}^{-1}\, v\, \Big(\, v\,A_{t,\tau}^{-1}\, v -
      v\, \tilA_{t,\tau}^{-1}\, v\, \Big)
   \ + \
   t^{2k}\, \Big(\, v\, A_{t,\tau,v}^{-1}\, v -
       v\, \tilA_{t,\tau,v}^{-1}\, v\, \Big)\,  v\, \tilA_{t,\tau}^{-1}\, v
   \\ + \
    \Big(\, v\, A_{t,\tau,v}^{-1}\, v -
             v\, \tilA_{t,\tau,v}^{-1}\, v\, \Big)
\end{multline}
As in \refe{Av-tilAv}, \refl{asexp+v} implies that for all $N\in\NN$ the traces of the
second and the third summands in the right hand side of \refe{vA-Avv} behave as
$o(t^{-N})$ when $t\to\infty$ uniformly in $\tau\in[0,1]$. Thus
\eq{TrvA-Avv}
    \Big|\, \Tr\, \Big[\, v\, A_{t,\tau}^{-1}\, v -
             v\, \tilA_{t,\tau}^{-1}\, v\, \Big]\, \Big|
    \ \le \
   t^{2k}\, \big\|\, v\, A_{t,\tau,v}^{-1}\, v\, \big\|\cdot
    \Big|\, \Tr\, \Big[\, v\,A_{t,\tau}^{-1}\, v -
      v\, \tilA_{t,\tau}^{-1}\, v\, \Big]\, \Big| \ + \ o(t^{-N}).
\end{equation}
Clearly, \/
 \(\displaystyle
    t^{2k}\, \big\|\, v\, A_{t,\tau,v}^{-1}\, v\, \big\|
    \ \le \ \frac{\eps^2t^{2k}}{\eps^2t^{2k}+|t|^{-k}}
    \ \le \  1- \frac{\eps^{-2}|t|^{-3k}}2.
 \) \/
Hence, from \refe{TrvA-Avv} we conclude
\eq{TrvA-Avv2}
    \Big|\, \Tr\, \Big[\, v\, A_{t,\tau}^{-1}\, v -
             v\, \tilA_{t,\tau}^{-1}\, v\, \Big]\, \Big|
    \ \le \ o(t^{3k-N}).
\end{equation}
Taking $N>5k$ we obtain from \refe{I=} and \refe{TrvA-Avv2} that $I_{t,\tau,v}=
o(t^{-2k})$ uniformly in $\tau\in [0,1]$.
\eprf


\section{Proof of Cheeger-M\"uller theorem}\Label{S:prRS}

Recall from \refc{RS=}, that to prove \reft{RS} it is enough to show that $R(M,F,f)=0$.
We will use the notation of \refs{witten}. Clearly $\lapt= \Del_{-f,-t}$. Hence,
\eq{rho-}
    \rho^{RS}_{la}(f,-t) \ = \ \rho^{RS}_{la}(-f,t).
\end{equation}
Recall that the number $R(M,F,f)$ is defined in \refe{R}. Let $\tilM,\tilF,\tilf$ be as
in \refss{2mflds}. It follows from \refc{RS=} that $R(M,F,f)- R(\tilM,\tilF,\tilf)$ is
equal to the free term of the asymptotic expansion of $\log\rsttl-
\log{}\rho^{RS}_{la}(\tilf,t)$. Hence, from \reft{compare} and \refe{rho-}, we conclude
\[
   \big[\, R(M,F,f)- R(\tilM,\tilF,\tilf)\, \big] \ + \
    \big[\, R(M,F,-f)- R(\tilM,\tilF,-\tilf)\, \big] \ = \ 0.
\]
Since $R(M,F,f)$ is independent of $f$, cf. \refe{R}, we obtain
\eq{compareR}
    R(M,F,f) \ = \ R(\tilM,\tilF,\tilf).
\end{equation}

\lem{product}
Suppose $N$ is a compact manifolds of even dimension. Let $\oF$ be the flat Hermitian
vector bundle induced by $F$ on the product $M\times{N}$. Fix a generalized triangulation
$(g^{TN},f^N)$ on $N$, cf. \refd{crpoint}. Let $\of$ be the Morse function on
$M\times{N}$ defined by the formula $\of(x,y)= f(x)+f^N(y)$, where $x\in M, y\in N$. Then
\eq{product}
    \log\rho^{RS}_{la}(\of,t) \ = \
    \chi(N)\log\rho^{RS}_{la}(f,t),
\end{equation}
where $\chi(N)$ is the Euler characteristic of $N$.
\elem
The proof is a verbatim repetition of the proof of Theorem~2.5 in \cite{RaSi1} and will
be omitted (In \cite{RaSi1}, the equality \refe{product} is proven with $\rho^{RS}_{la}$
replaced by the ``full" Ray-Singer torsion $\rho^{RS}$). Using \refc{RS=} and
\refl{product}, we obtain
\eq{product2}
    R(M\times N,\oF,\of) \ = \ \chi(N)\, R(M,F,f).
\end{equation}

Substituting in \refe{product2},  $N=S^2$ and $N=S^1\times{S^1}$ we obtain respectively
\eq{product3}
    R(M\times S^2,\oF,\of) \ = \ 2R(M,F,f); \qquad
    R(M\times S^1\times S^1,\oF,\of) \ = \ 0.
\end{equation}

Using the results of \cite[\S5]{Milnor-hcob} (see also Lemma~4.2 of \cite{BFK3}), we see
that there exist generalized triangulations $(g^{M\times{S^2}},f_1)$ on $M\times{S^2}$
and $(g^{M\times{S^1}\times{S^1}},f_2)$ on $M\times{S^1}\times{S^1}$, such that the
functions $f_1$ and $f_2$ have the same critical points structure, cf. \refd{samecp}.
Hence, by \refe{compareR}, we have
\eq{S1=S2}
    R(M\times S^2,\oF,f_1) \ = \ R(M\times S^1\times S^1,\oF,f_2).
\end{equation}
From \refe{product3}, \refe{S1=S2} and the fact that $R$ is independent of the choice of
the Morse function, we obtain $R(M,F,f)=0$. \hfill$\square$

\bibliographystyle{plain}
\providecommand{\bysame}{\leavevmode\hbox to3em{\hrulefill}\thinspace}

\end{document}